\documentclass[a4paper,12pt]{amsart}
\usepackage{amsmath}
\usepackage{amsthm}
\usepackage{amscd}
\usepackage{amssymb}
\usepackage{amsfonts}
\usepackage{latexsym}
\usepackage[mathscr]{eucal}

\newtheorem{thm}{Theorem}[section]
\newtheorem{prop}[thm]{Proposition}
\newtheorem{lemma}[thm]{Lemma}

\def\graph{\mathop{\rm {graph}}\nolimits}

\def\Im{\mathop{\rm {Im}}\nolimits}

\def\Index{\mathop{\rm {Index}}\nolimits}
\def\Ker{\mathop{\rm {Ker}}\nolimits}

\begin{document}

\title{exotic indecomposable systems of four subspaces
in a Hilbert space}
\author{Masatoshi Enomoto}
\address[Masatoshi Enomoto]{College of Business Administration 
and Information Science, 
Koshien University, Takarazuka, Hyogo 665, Japan}      

\author{Yasuo Watatani}
\address[Yasuo Watatani]{Department of Mathematical Sciences, 
Kyushu University, Hakozaki, 
Fukuoka, 812-8581,  Japan}

\dedicatory{Dedicated to Professor Masahiro Nakamura on his 
88th birthday} 
\maketitle
\begin{abstract}We study the relative position of four 
subspaces in a Hilbert space. 
For any positive integer $n$, 
we give an example of exotic indecomposable system ${\mathcal S}$ 
of four subspaces in a Hilbert space whose defect is 
$ \frac{2n+1}{3}$.
By an exotic system, we mean a system which is not 
isomorphic to any closed operator system under any 
permutation of subspaces. We construct the examples
by a help of certain nice sequences used by 
Jiang and Wang in their study of strongly irreducible 
operators.


\medskip\par\noindent
AMS SUBJECT CLASSIFICATION: 46C07, 47A15, 15A21, 16G20, 16G60.

\end{abstract}

\section{Introduction} 

Many problems of linear algebra 
can be reduced to the classification of the systems of $n$ subspaces 
in a finite-dimensional vector space. Nazarova \cite{N} and 
 Gelfand-Ponomarev \cite{GP} completely classified 
indecomposable systems of four subspaces in a finite dimensional 
vector space. On the other hand, in operator theory, 
Halmos initiated the study of transitive lattices of subspaces, 
see for example \cite{Ha}. Transitive lattices give 
transitive systems of subspaces.  Transitive system of subspaces 
in a finite dimensional space had been studied by Brenner in \cite{B}. 

In \cite{EW} we started to investigate systems of $n$ 
subspaces in an infinite dimensional Hilbert space considering 
an analogy with subfactor theory invented by Jones \cite{J}. 
As a building block,  we investigate indecomposable systems 
of $n$ subspaces in the sense that the system can not be isomorphic 
to a direct sum of two  non-zero systems. Recently Moskaleva and 
Samoilenko \cite{MS} study a relation between systems of $n$-subspaces and 
representations of *-algebras generated by projections. 

Let $H$ be a Hilbert space and $E_1, \dots E_n$ $n$ subspaces 
in $H$.  Then we say that  ${\mathcal S} = (H;E_1, \dots , E_n)$  
is a system of $n$ subspaces in $H$ or a $n$-subspace system in $H$.
A system ${\mathcal S}$ is called {\it indecomposable} 
if ${\mathcal S}$ is not be decomposed into a nontrivial direct sum.
  
For any bounded linear operator $A$ on a Hilbert space $K$, we 
associate an operator system ${\mathcal S}_A$ of four subspaces in 
$H = K \oplus K$ by 
\[
{\mathcal S}_A = (H;K\oplus 0,0\oplus K,\graph A, \{(x,x) ; x \in K\}).
\]
Two such operator systems ${\mathcal S}_A$ and ${\mathcal S}_B$  
are isomorphic 
if and only if the two operators $A$ and $B$ are similar. The 
direct sum of operator systems corresponds to the direct sum of 
the operators. In this sense the study of operators is 
included into the study of relative positions of four subspaces. 
In particular in a finite dimensional space, Jordan blocks correspond 
to indecomposable systems. Moreover in an infinite dimensional 
Hilbert space, an operator system ${\mathcal S}_A$ is indecomposable 
if and only if $A$ is strongly irreducible.  Recall that 
an operator $A\in B(K)$ is called strongly irreducible if 
there are no non-trivial invariant subspaces $M$ and $N$ 
of $A$ such that $M \cap N = 0$ and $M + N = K$. 
A strongly irreducible operator is an 
infinite-dimensional analog of a Jordan block. We refer
a good monograph \cite{JW} by Jiang and Wang on 
strongly irreducible operators. 

In \cite{EW} we discovered some examples of exotic indecomposable systems 
$\mathcal {S}$ 
of four subspaces in a Hilbert space.  
By an exotic system, we mean a system which is not 
isomorphic to any closed operator system ${\mathcal S}_A$ under any 
permutation of subspaces. 

  Gelfand and Ponomarev introduced an integer 
valued  invariant $\rho ({\mathcal S})$, 
called  
{\it defect},  for a system 
${\mathcal S} = (H ; E_1, E_2, E_3, E_4)$ of four subspaces by 
\[
\rho ({\mathcal S}) = \sum _{i=1} ^4 \dim E_i - 2\dim H.  
\]
They showed that if ${\mathcal S}$ is indecomposale, then 
the defect $\rho ({\mathcal S})$ is one of $\{-2,-1,0,1,2\}$. 

 We extended the notion of defect to a certain class of 
systems of four subspaces in  an infinite dimensional Hilbert space 
using  Fredholm index in \cite{EW}. We showed that the defect 
for indecomposable systems of 
four subspaces takes any value in  $\mathbb Z/3$. 
These values are attained by bounded operator systems. 
In fact the exotic systems constructed 
in \cite{EW} have the defect $\rho ({\mathcal S}) = 1$. 

The aim of the paper is to give new examples 
 of {\it exotic} indecomposable systems $\mathcal {S}$ 
of four subspaces in a Hilbert space with the defect 
$\rho ({\mathcal S}) = \frac{2n+1}{3}$ 
for any positive integer $n$.  
We construct these examples
by a help of certain nice sequences used by 
Jiang and Wang in their study of  strongly irreducible 
operators in \cite{JW}.

\section{relative position of subspaces}

We study the relative position of $n$ subspaces
in a separable  Hilbert space. 
Firstly we recall some basic facts in \cite{EW}. 
Let $H$ be a Hilbert space and $E_1, \dots ,E_n$ be $n$ subspaces 
in $H$.  Then we say that  ${\mathcal S} = (H;E_1, \dots , E_n)$  
is a system of $n$-subspaces in $H$ or an $n$-subspace system in $H$. 
Let ${\mathcal T} = (K;F_1, \dots , F_n)$  
be  another system of $n$-subspaces in a Hilbert space $K$. Then  
$\varphi : {\mathcal S} \rightarrow {\mathcal T}$ is called a 
homomorphism if $\varphi : H \rightarrow K$ is a bounded linear 
operator satisfying that  
$\varphi(E_i) \subset F_i$ for $i = 1,\dots ,n$. And 
$\varphi : {\mathcal S} \rightarrow {\mathcal T}$
is called an isomorphism if $\varphi : H \rightarrow K$ is 
an invertible (i.e., bounded  bijective) linear 
operator satisfying that  
$\varphi(E_i) = F_i$ for $i = 1,\dots ,n$. 
We say that systems ${\mathcal S}$ and ${\mathcal T}$ are 
{\it isomorphic} if there is an isomorphism  
$\varphi : {\mathcal S} \rightarrow {\mathcal T}$. This means 
that the relative positions of $n$ subspaces $(E_1, \dots , E_n)$ in $H$ 
and   $(F_1, \dots , F_n)$ in $K$ are same under disregarding angles. 
We say that systems ${\mathcal S}$ and ${\mathcal T}$ are 
{\it unitarily equivalent } if the above isomorphism 
$\varphi : H \rightarrow K$ can be chosen to be a unitary. 
This means that the relative positions of 
$n$ subspaces $(E_1, \dots , E_n)$ in $H$ 
and  $(F_1, \dots , F_n)$ in $K$ are same with preserving the angles
between the subspaces.

We denote by 
$Hom(\mathcal S, \mathcal T)$ the set of homomorphisms of 
$\mathcal S$ to $\mathcal T$ and  
$End(\mathcal S) := Hom(\mathcal S, \mathcal S)$ 
the set of endomorphisms on $\mathcal S$. 

For two \ systems $\mathcal S=(H;E_{1},...,E_{n})$ and
$\mathcal T=(K;F_{1},...,F_{n})$ of $n$ subspaces in $H$, 
their direct sum $\mathcal S\oplus \mathcal T$ is
defined by
$$
\mathcal S\oplus \mathcal T:=(H\oplus K;E_{1}\oplus F_{1},...,E_{n}\oplus F_{n}).$$

\noindent
{\bf Definition.} 
A system $\mathcal S =(H;E_{1},...,E_{n})$ of $n$ subspaces is called 
{\it decomposable} if the
systems $\mathcal S$ is isomorphic to \ a direct sum of two \ non-zero \ systems. A
system $\mathcal S =(H;E_{1},...,E_{n})$ is said to \ be 
{\it indecomposable}  if it is not decomposable. A system 
$\mathcal S$ is indecomposable if and only if 
$Idem (\mathcal S) := \{V \in End (\mathcal S) ; V^2 = V\} = \{0,I\}$. 
A system 
$\mathcal S$ is said to be {\it transitive}  if 
$End (\mathcal S) = {\mathbb C}I$. 

\smallskip

Transitive systems in a finite dimensional space were studied by S. Brenner \cite{B}. On the other hand,  Halmos \cite{Ha} initiated the study of transitive 
lattices of subspaces in Hilbert spaces, which give transitive systems. 
Some interesting examples were obtained by 
Harrison-Radjavi-Rosenthal \cite{HRR} and Hadwin-Longstaff-Rosenthal \cite{HLR}.
We have a close relation between systems of subspaces and operators. In fact   
we can associate a system of four subspaces for any operator. 

\medskip
\noindent
{\bf Definition.} 
We say that a system $\mathcal S=(H;E_{1},...,E_{4})$ of four subspaces is
a {\it closed operator system}  if there exist Hilbert spaces 
$K_{1},$ $K_{2}$ and closed
operators $T: K_{1}\supset D(T)\rightarrow K_{2}$,  $S: K_{2}\supset
D(S)\rightarrow K_{1}$  
such that $H=K_{1}\oplus$ $K_{2}$, 
$E_{1}=K_{1}\oplus 0$, $E_{2}=0\oplus K_{2}$, $E_{3}=\{(x,Tx);x\in
D(T)\}$ and $E_{4}=\{(Sy,y);y\in D(S)\}$. Here $D(T)$ is the domain of $T$.
In particular,  if $T$ and $S$ are bounded operators with
 $D(T) = K_1$ and $D(S) = K_2$, then we say that 
$\mathcal S=(H;E_{1},...,E_{4})$ is a bounded operator system. 
We denote it by ${\mathcal S}_{T,S}$. We put 
${\mathcal S}_{T} := {\mathcal S}_{T,I}$ and call it a bounded operator 
system associated with a single operator $T$. 
Two such operator systems ${\mathcal S}_A$ and ${\mathcal S}_B$  
are isomorphic 
if and only if the two operators $A$ and $B$ are similar. 
Moreover in an infinite dimensional 
Hilbert space, a bounded operator system ${\mathcal S}_A$ is indecomposable 
if and only if $A$ is strongly irreducible.  

\medskip
\noindent
{\bf Definition.} Let ${\mathcal S} = (H ; E_1, E_2, E_3, E_4)$ 
be a system of four subspaces. For any distinct $i,j = 1,2,3,4$, 
define an adding operator 
\[
A_{ij} : E_i \oplus E_j \ni (x,y) \rightarrow x+y \in H. 
\] 
Then 
\[
\Ker A_{ij} = \{(x,-x) \in E_i \oplus E_j ; x \in E_i \cap E_j \}
\]
and 
\[
\Im A_{ij} = E_i + E_j .
\]
We say ${\mathcal S} = (H ; E_1, E_2, E_3, E_4)$ is a
{\it Fredholm} system if  $A_{ij}$ is a Fredholm operator for any 
$i,j = 1,2,3,4$ with $i \not= j$. Then $\Im A_{ij} = E_i + E_j$ is closed 
and 
\[
\Index A_{ij} = \dim \Ker A_{ij}  - \dim \Ker A_{ij}^* 
= \dim (E_i \cap E_j) - \dim ((E_i + E_j)^{\perp}). 
\]

\medskip
\noindent
{\bf Definition.} We say 
${\mathcal S} = (H ; E_1, E_2, E_3, E_4)$ is a
{\it quasi-Fredholm} system if $E_i \cap E_j$ and $(E_i + E_j)^{\perp}$ 
are finite-dimensional for any $i \not= j$. In the case we define 
the {\it defect} $\rho ({\mathcal S})$ of ${\mathcal S}$ by 

$$
\rho ({\mathcal S}) := 
\frac{1}{3} \sum _{1\leq i<j\leq 4} 
(\dim (E_i \cap E_j) - \dim (E_i + E_j)^{\perp})
$$
which coincides with the Gelfand-Ponomarev original defect if 
$H$ is finite-dimensional. 
Moreover, if ${\mathcal S}$ is a Fredholm system, then 
it is a quasi-Fredholm system and 
\[
\rho ({\mathcal S}) = \frac{1}{3} \sum _{1\leq i<j\leq 4} \Index A_{ij} . 
\]

\section{construction of examples}

Consider a Hilbert space  $L=\ell^{2}({\mathbb N})$. 
Let $\{e_1, e_2, e_3, \dots \}$ be  a canonical basis . 
For a bounded sequence $w = (w(n))_n$, we define a 
backward weighted shift $B_w \in B(\ell^{2}({\mathbb N}))$ 
of weight $w$ by 
$$
B_{w}e_n = w(n-1)e_{n-1}, \ (n \geq 2)  \ \text{ and } \ 
B_{w}e_1 = 0. 
$$
Thus for $x = (x(n))_n \in \ell^{2}({\mathbb N})$, we have 
$(B_wx)(n) = w(n)x(n+1)$ for $n = 1,2,3,\dots $.

We borrow a family of sequences $a_1=(a_1(n))_n, a_2 =(a_2(n))_n, 
a_3=(a_3(n))_n,\dots $  used by Jiang and Wang in \cite[p.93-94]{JW} 
as follows: 

Define a sequence $c = (c(n))_n$ of positive numbers 
and an increasing sequence $n_1 < n_2 <n_3 < \dots   $
of natural numbers as follows:

Put $c(1) = 2 = \frac{1+1}{1} >1 $ and $n_1 = 1$. 
There exists $n_2 \in {\mathbb N}$ with $n_1 < n_2$ such that 
$$
\frac{1+1}{1} \prod_{k=n_1+1}^{n_2} \frac{k}{k+1} < \frac{1}{2}.
$$
Put $c(k) = \frac{k}{k+1}$ for $k = n_1 +1 = 2, \dots, n_2$. 
There exists $n_3 \in {\mathbb N}$ with $n_2 < n_3$ such that 
$$
\frac{1+1}{1} \prod_{k=n_1+1}^{n_2}\frac{k}{k+1} 
\prod_{k=n_2+1}^{n_3}\frac{k+1}{k} > 3 . 
$$
Put $c(k) = \frac{k+1}{k}$ for $k = n_2 +1, \dots, n_3$. 
We continue in this fashion to obtain an increasing sequence 
$n_1 < n_2 <n_3 < \dots   $
of natural numbers and a sequence $c = (c(n))_n$ of 
positive numbers such that 
$$
c(k) = 
\begin{cases}
\frac{k+1}{k}, \ \ (k = n_1 = 1,  n_2 +1 \leq k \leq n_3, 
                 n_4 +1 \leq k \leq n_5, \dots   ) \\
\frac{k}{k+1}, \ \ (n_1 + 1 \leq k \leq n_2,  
                 n_3 +1 \leq k \leq n_4, 
                 n_5 +1 \leq k \leq n_6, \dots   ),  
\end{cases}
$$
and 
$$
\prod_{k=1}^{n_j} c(k) 
\begin{cases}
> j \ \ (j \text{ is odd }) \\
< \frac{1}{j} \ \ (j  \text{ is even. })
\end{cases}
$$
Then $\frac{2}{3} \leq c(k) \leq 2.$ 
%Put 
%$$
%N_1 = \{ k \in {\mathbb N} \ ;\ 1 \leq k \leq n_1, 
%        n_2 +1 \leq k \leq n_3, n_4 +1 \leq k \leq n_5, \dots \ \}
%$$
%$$
%N_2 =  \{ k \in {\mathbb N} \ ; \  n_1 +1 \leq k \leq n_2, \
%        n_3 +1 \leq k \leq n_4, \\ 
%   n_5 +1 \leq k \leq n_6, \dots \}
%$$
Define $a_1(k) \equiv 1$, $a_2(k) = c(k)^{\frac{1}{2}}$ and 
$$
a_i(k) = c(k)^{ \frac{1}{2} + \frac{1}{4}+ \dots + \frac{1}{2^{i-1}}} 
       = c(k)^{1-2^{1-i}} \text{ for } i \geq 2
$$
Then we have the following lemma:

\begin{lemma}(Jiang and Wang \cite{JW}) There exists 
a family of sequences $a_1=(a_1(n))_n, a_2 =(a_2(n))_n, 
a_3=(a_3(n))_n,\dots $, of positive numbers satisfying 

\begin{enumerate}
\item $\frac{2}{3} \leq a_i(k) \leq 2,$  
\item $\lim_{k \rightarrow \infty} a_i(k) = 1$, \
$\lim_{n \rightarrow \infty} \prod_{k=1}^{n} 2a_i(k) = \infty ,$
\item $\limsup_{n \rightarrow \infty} \prod_{k=1}^{n} 
\frac{a_i(k)}{a_j(k)} = \infty,  \ (i \not= j)$, 
\item $\liminf_{n \rightarrow \infty} \prod_{k=1}^{n} 
\frac{a_i(k)}{a_j(k)} = 0,  \ (i \not= j)$.
\item the point spectrum $\sigma _p(B_{a_i})$ contains
 $\{\lambda \in {\mathbb C} : | \lambda | <1 \}$. 
\end{enumerate}
\label{lemma;sequences}
\end{lemma}

We shall construct our examples.  We fix 
a family of sequences $a_1=(a_1(n))_n, a_2 =(a_2(n))_n, 
a_3=(a_3(n))_n,\dots $ of positive numbers defined in 
the above  Lemma \ref{lemma;sequences}. Put $w_k = 2a_k$ for  
$k = 1,2,\dots $. 
Consider a sequence of backward weighted shifts 
$B_{w_1}, B_{w_2},B_{w_3},...$  on $L=\ell^{2}({\mathbb N})$.
Let $S \in B(L)$ be a unilateral shift.  For a fixed natural 
number $N$,  define  $K = L \oplus \dots \oplus L$ 
($N + 1$ times) and $H = K \oplus K$.  In the below we sometimes 
use symbol $(x\oplus y) \in K \oplus K$ instead of 
$(x,y) \in K \oplus K$ for the sake of convenience of notation.  
We consider an operator 
$$
T=\left(
\begin{array}
[c]{ccccc}%
B_{w_{1}} & I & O & ... & O\\
O & B_{w_{2}} & I & \ddots & \vdots\\
O & O & \ddots & \ddots & O\\
\vdots & \ddots & \ddots & B_{w_{N}} & I\\
O & O & \cdots & O & S
\end{array}
\right)
\in B(K)  
$$

Let $E_1=K \oplus 0$, 
$E_2= 0 \oplus  K$, $E_{4}=\{x\oplus x \in K \oplus K ; x \in K\}$ and 
$$
E_3=\{ x \oplus Tx \in K \oplus K ;  x \in K \} 
+ {\mathbb C}((0,\dots,0)\oplus (0,\dots,0,e_1)). 
$$ 

Consider a system $\mathcal{S}_{w,N} =(H;E_1,E_2,E_3,E_4)$. 
We shall show that $\mathcal{S}_{w,N}$ is indecomposable and is not
isomorphic to any closed operator systems under any permutation.  
We could  regard that the system  $\mathcal{S}_{w,N}$ is a 
one-dimensional \lq\lq deformation" of an operator system, since 
$E_3 = \graph T +  {\mathbb C}((0,\dots,0)\oplus (0,\dots,0,e_1)).$

\begin{thm}
The above system $\mathcal{S}_{w,N}$ of four subspaces is indecomposable. 
\label{thm:exotic examples} 
\end{thm}
\begin{proof}In order to make the notation simple, we shall prove 
the theorem in case $N = 3$. The general $N$ case will be 
proved similarly. 
Let $V \in End (\mathcal{S}_{w,N})$ satisfy $V^2 = V$. It is enough 
to show that $V = O$ or $V = I$ for $\mathcal{S}_{w,N}$ to be 
indecomposable. 
Since $V(E_i) \subset E_i$ for $i = 1,2,4$, we have 
\[
V = 
    \begin{pmatrix}
     A & O \\
     O & A \\ 
    \end{pmatrix} 
  \in B(H) \ \ \text{ for some } \ A \in B(K).
\]
It is sufficient to prove that $A = O$ or $A = I$.  
We may write $A = (A_{ij})_{ij}$ as operator matrix , 
where $A_{ij} \in B(L)$ and 
$i, j = 1,2,3,4.$  
Thus we have 
\[
T = 
\begin{pmatrix}
     B_{w_1} & I & O & O \\
     O & B_{w_2} & I & O \\
     O & O & B_{w_3} & I \\
     O & O & O & S
\end{pmatrix}
\quad \text{and }
A = 
\begin{pmatrix}
     A_{11} & A_{12} & A_{13} & A_{14} \\
     A_{21} & A_{22} & A_{23} & A_{24} \\
     A_{31} & A_{32} & A_{33} & A_{34} \\
     A_{41} & A_{42} & A_{43} & A_{44}
\end{pmatrix}
.
\]
Since $E_3 = \graph T +  {\mathbb C}((0,0,0,0)\oplus (0,0,0,e_1) )$, 
$E_3$ is spanned by 
\[
\left\{
\begin{array}
[c]{c}%
(e_{1},0,0,0)\oplus(0,0,0,0), 
\ (e_{n},0,0,0)\oplus(w_{1}(n-1)e_{n-1},0,0,0),\\
(0,e_1,0,0)\oplus(e_1,0,0,0), 
\ 
(0,e_{n},0,0)\oplus(e_{n},w_{2}(n-1)e_{n-1},0,0),\\
(0,0,e_1,0)\oplus(0,e_1,0,0), 
\ (0,0,e_{n},0)\oplus(0,e_{n},w_{3}(n-1)e_{n-1},0),\\
(0,0,0,e_{k})\oplus(0,0,e_{k},e_{k+1}),\\
(0,0,0,0)\oplus(0,0,0,e_{1}), 
\end{array}
;  
\begin{array}
[c]{c}%
n \geq 2, \\
 k \geq 1 
\end{array}
\right\}. 
\]
We may write 
$$
E_3 = 
\left\{ \left(
\begin{array}
[c]{c}%
(x_{1}(n))_n\\
(x_{2}(n))_n\\
(x_{3}(n))_n\\
(x_{4}(n))_n
\end{array}
\right)  \oplus \left(
\begin{array}
[c]{c}%
(x_{1}(n+1)w_{1}(n)+x_{2}(n))_n\\
(x_{2}(n+1)w_{2}(n)+x_{3}(n))_n\\
(x_{3}(n+1)w_{3}(n)+x_{4}(n))_n\\
(y, (x_4(n))_n)
\end{array}
\right)  
; 
\begin{array}
[c]{c}%
x_1,x_2 \in \ell^{2}({\mathbb N}) \\
x_3,x_4 \in \ell^{2}({\mathbb N}) \\
y \in {\mathbb C}
\end{array}
\right\}
$$
We need several lemmas in the below  to complete the proof.
\end{proof}

\begin{lemma}Let $P\in B(\ell^{2}({\mathbb N}))$ be an operator 
of the form $P = \lambda I + N$ for some 
$\lambda \in {\mathbb C}$ and  an upper (or lower)
triangular matrix $N \in B(\ell^{2}({\mathbb N}))$ 
with zero diagonal.  Assume that $P$ is an idempotent, 
then $P = O$ or $P = I$. 
\label{lemma;known}  
\end{lemma}
\begin{proof}
This is a known fact. See for example Lemma 10.1 in \cite{EW}.
\end{proof}

\begin{lemma} We have that 
$A_{41}(k,n) = 0$ for any $k,n \geq 1$,
 $A_{31}(k,n) = 0$ for any $k \geq n+1$, 
$A_{21}(k,n) = 0$ for any $k \geq n+2$, 
and $A_{11}(k,n) = 0$ for any $k \geq n+3$. 
In particular $A_{41} = O$. 
\label{lemma;triangular}
\end{lemma}
\begin{proof}
Since $u = (e_{1},0,0,0)\oplus(0,0,0,0) \in E_3$, we have 
\begin{align*}
Vu & =  
\left(
\begin{array}
[c]{c}%
A_{11}e_1\\
A_{21}e_1\\
A_{31}e_1\\
A_{41}e_1
\end{array}
\right)  \oplus \left(
\begin{array}
[c]{c}%
0\\
0\\
0\\
0
\end{array}
\right)  \\
& = 
\left(
\begin{array}
[c]{c}%
(x_{1}(k))_k\\
(x_{2}(k))_k\\
(x_{3}(k))_k\\
(x_{4}(k))_k
\end{array}
\right)  \oplus \left(
\begin{array}
[c]{c}%
(x_{1}(k+1)w_{1}(k)+x_{2}(k))_k\\
(x_{2}(k+1)w_{2}(k)+x_{3}(k))_k\\
(x_{3}(k+1)w_{3}(k)+x_{4}(k))_k\\
(y, (x_4(k))_k)
\end{array}
\right)  
\in E_3.
\end {align*}
for some 
$x_1,x_2,x_3,x_4 \in \ell^{2}({\mathbb N})$ and 
$y \in {\mathbb C}$. 
Then $x_4(k) = 0$ for $k \geq 1$. Thus $A_{41}(k,1) = (A_{41}e_1)(k) 
= x_4(k) = 0$.  

Since 
$x_3(k+1)w_3(k) = x_3(k+1)w_3(k) + x_4(k) = 0$ and $w_3(k) > 0$, 
we have $x_3(k+1) = 0$ for $k \geq 1$, i.e., 
$A_{31}(k,1) = x_3(k) = 0$ for $k \geq 2$.  

Since 
$x_2(k+1)w_2(k) = x_2(k+1)w_2(k) + x_3(k) = 0 $ 
for $k \geq 2$ and $w_2(k) > 0$, 
we have $x_2(k+1) = 0$ for $k \geq 2$, i.e., 
$A_{21}(k,1) = x_2(k) = 0$ for $k \geq 3$.  

Since 
$x_1(k+1)w_1(k) = x_1(k+1)w_1(k) + x_2(k) = 0 $ 
for $k \geq 3$ and $w_1(k) > 0$, 
we have $x_1(k+1) = 0$ for $k \geq 3$, i.e., 
$A_{11}(k,1) = x_1(k) = 0$ for $k \geq 4$.  
Thus the statement of the lemma is proved for  $n = 1$. 

Moreover, $x_2(2)w_2(1) + x_3(1) = 0$ implies that 
$A_{21}(2,1) w_2(1) + A_{31}(1,1) = 0$. \\
And $x_1(2)w_1(1) + x_2(1) = 0$ implies that 
$A_{11}(2,1) w_1(1) + A_{21}(1,1) = 0$. 
And 
$x_1(3)w_1(2) + x_2(2) = 0$ implies that 
$A_{11}(3,1) w_1(2) + A_{21}(2,1) = 0$. 

We shall prove the lemma by induction on $n$. 
Assume that the statement of the Lemma holds for the $n$-th 
column of $A_{11}, A_{21}, A_{31}, A_{41}$ . 
We shall prove it for $n+1$.  

Since $u = (e_{n+1},0,0,0)\oplus(w_1(n)e_n,0,0,0) \in E_3$, we have 
\begin{align*}
Vu & =  
\left(
\begin{array}
[c]{c}%
A_{11}e_{n+1}\\
A_{21}e_{n+1}\\
A_{31}e_{n+1}\\
A_{41}e_{n+1}
\end{array}
\right)  \oplus \left(
\begin{array}
[c]{c}%
w_1(n)A_{11}e_{n}\\
w_1(n)A_{21}e_{n}\\
w_1(n)A_{31}e_{n}\\
w_1(n)A_{41}e_{n}
\end{array}
\right)  \\
& = 
\left(
\begin{array}
[c]{c}%
(x_{1}(k))_k\\
(x_{2}(k))_k\\
(x_{3}(k))_k\\
(x_{4}(k))_k
\end{array}
\right)  \oplus \left(
\begin{array}
[c]{c}%
(x_{1}(k+1)w_{1}(k)+x_{2}(k))_k\\
(x_{2}(k+1)w_{2}(k)+x_{3}(k))_k\\
(x_{3}(k+1)w_{3}(k)+x_{4}(k))_k\\
(y, (x_4(k))_k)
\end{array}
\right)  
\in E_3.
\end{align*}
for some 
$x_1,x_2,x_3,x_4 \in \ell^{2}({\mathbb N})$ and 
$y \in {\mathbb C}$. 
Since $(A_{41}e_n)(k) = A_{41}(k,n) = 0$ for any $k$ 
by the assumption of induction, 
$(y, (x_4(k))_k) = w_1(n)A_{41}e_{n} = 0$. Then    
$ A_{41}(k,n+1) = (A_{41}e_{n+1})(k) = x_4(k) = 0$. 

Since $(A_{31}e_n)(k) = A_{31}(k,n) = 0$ for any $k \geq n+1$ 
by the assumption of induction, 
$$
x_{3}(k+1)w_{3}(k) = x_{3}(k+1)w_{3}(k)+x_{4}(k) 
= w_1(n)(A_{31}e_{n})(k) = 0. 
$$ 
Because 
$w_3(k) > 0$, 
we have $x_3(k+1) = 0$ for $k \geq n+1$, i.e., 
$$
A_{31}(k,n) = (A_{31}e_n)(k) =  x_3(k) = 0
\text{ for } k \geq (n+1) +1 .
$$
Since $A_{21}(k,n) = 0$ for any $k \geq n+2$ 
by the assumption of induction, 
$$
x_{2}(k+1)w_{2}(k) = x_{2}(k+1)w_{2}(k)+x_{3}(k) 
= w_1(n)(A_{21}e_{n})(k) = 0. 
$$ 
Because 
$w_2(k) > 0$, 
we have $x_2(k+1) = 0$ for $k \geq n+2$, i.e., 
$$
A_{21}(k,n) = (A_{21}e_n)(k) =  x_2(k) = 0 
\text{ for } k \geq (n+1) +2. 
$$
Since $A_{11}(k,n) = 0$ for any $k \geq n+3$ 
by the assumption of induction, 
$$
x_{1}(k+1)w_{1}(k) = x_{1}(k+1)w_{1}(k)+x_{2}(k) 
= w_1(n)(A_{11}e_{n})(k) = 0. 
$$ 
Because 
$w_1(k) > 0$, 
we have $x_1(k+1) = 0$ for $k \geq n+3$, i.e., 
$$
A_{11}(k,n) = (A_{11}e_n)(k) =  x_1(k) = 0
\text{ for } k \geq (n+1) +2.
$$ 

This finishes the proof by induction. 
\end{proof}

\begin{lemma} $A_{31} = O$ and $A_{21} = O$. 
\end{lemma}
\begin{proof}
From  the proof in Lemma \ref{lemma;triangular}, 
$A_{31}(k,n+1) = x_3(k) $ and 
$$
w_1(n)A_{31}(k,n) = x_3(k+1)w_3(k) + x_4(k) = x_3(k+1)w_3(k). 
$$
Hence 
$$
A_{31}(k+1,n+1) = \frac{w_1(n)}{w_3(k)} A_{31}(k,n) 
\text{ for any }  n, k. 
$$
Therefore for any $j \geq 1$
$$
A_{31}(1+n,j+n) 
= \prod_{k=1}^{n} \frac{w_1(j+k-1)}{w_3(1+k-1)} A_{31}(1,j).
$$
Recall that $\frac{4}{3}\leq w_1(k) \leq 4$, 
$\frac{4}{3}\leq w_3(k) \leq 4$ and 
 $\limsup_{n \rightarrow \infty} \prod_{k=1}^{n} 
\frac{w_1(k)}{w_3(k)} = \infty$ by Lemma \ref{lemma;sequences}.  
Since $\|A_{31}\| < \infty$, we have $A_{31}(1,j) = 0$. 
Furthermore  $A_{31}(1+n,j+n) = 0$ for any $j, n$, 
i.e., $A_{31}(k,n) = 0$ for any $k \leq n$. 
By Lemma \ref{lemma;triangular}
$A_{31}(k,n) = 0$ for any $k \geq n+1$. 
Therefore $A_{31} = O$.

Similarly we have 
$$
A_{21}(k+1,n+1) = \frac{w_1(n)}{w_2(k)} A_{21}(k,n) 
\text{ for any }  n, k. 
$$
By a similar argument we also have $A_{21} = O$.
\end{proof}

\begin{lemma} $A_{11} = O$ or $A_{11} = I$. 
\label{lemma;diagonal}
\end{lemma}
\begin{proof}
From the proof in Lemma \ref{lemma;triangular} and 
the additional fact that $x_2= A_{21}e_{n+1} = 0$, we have 
$A_{11}(k,n+1) = x_1(k) $ and 
$$
w_1(n)A_{11}(k,n) = x_1(k+1)w_1(k) + x_2(k) = x_1(k+1)w_1(k). 
$$
Hence 
$$
A_{11}(k+1,n+1) = \frac{w_1(n)}{w_1(k)} A_{11}(k,n) 
\text{ for any }  n, k. 
$$
Therefore for any $j \geq 1$
$$
A_{11}(j+n,1+n) 
= \prod_{k=1}^{n} \frac{w_1(1+k-1)}{w_1(j+k-1)} A_{11}(j,1).
$$
And we also have 
$$
A_{11}(1+n,1+n) = A_{11}(n,n).  
$$
Therefore the diagonal of $A_{11}$ is a constant, say $\lambda$. 
From the proof in Lemma \ref{lemma;triangular}, 
we have $A_{11}(3,1) = - \frac{A_{21}(2,1)}{w_1(2)} = 0$ and 
$A_{11}(2,1) = - \frac{A_{21}(1,1)}{w_1(1)} = 0$, 
because $A_{21} = O$.  Therefore  for any $n$ 
$$
A_{11}(n +2,n) 
= \prod_{k=1}^{n} \frac{w_1(1+k-1)}{w_1(3+k-1)} A_{11}(3,1) = 0.
$$ 
Similarly $A_{11}(n +1,n) = 0$. 
We also have $A_{11}(k,n) = 0$ for any $k \geq n+3$ by 
Lemma \ref{lemma;triangular}. 
Therefore $A_{11} = \lambda I + N$ for some 
$\lambda \in {\mathbb C}$ and  an upper 
triangular matrix $N \in B(\ell^{2}({\mathbb N}))$ 
with zero diagonal. Since $V$ is an idempotent, $A$ is an idempotent. 
Hence $A_{11}$ is also an idempotent, because  $A_{21} = A_{31} = A_{41} = O$. 
Thus $A_{11} = O$ or $A_{11} = I$ by Lemma \ref{lemma;known}. 
\end{proof}

In the below we shall show that if $A_{11} = O$(resp. $A_{11}= I$), 
then $V = O$ (resp. $V = I$). Replacing $V$ by $I-V$, it is enough to show 
that $A_{11} = O$ implies $V = O$ to prove Theorem 
\ref{thm:exotic examples}. 

\begin{lemma} Suppose that $A_{11} = O$. 
Then  $A_{42}(k,n) = 0$ for any $k,n \geq 1$,
 $A_{32}(k,n) = 0$ for any $k \geq n+1$, 
$A_{22}(k,n) = 0$ for any $k \geq n+2$, 
and $A_{12}(k,n) = 0$ for any $k \geq n+3$. 
In particular $A_{42} = O$. 
\end{lemma}
\begin{proof}
Since $u = (0,e_1,0,0)\oplus(e_1,0,0,0) \in E_3$ and 
the first column of $A$ is $0$, 
$$
Vu =  
\left(
\begin{array}
[c]{c}%
A_{12}e_1\\
A_{22}e_1\\
A_{32}e_1\\
A_{42}e_1
\end{array}
\right)  \oplus \left(
\begin{array}
[c]{c}%
0\\
0\\
0\\
0
\end{array}
\right)  \in E_3
$$
Since $u = (0,e_{n+1},0,0)\oplus(e_{n+1},w_2(n)e_n,0,0) \in E_3$ and 
the first column of $A$  is  $0$, 
we have 
$$
Vu =  
\left(
\begin{array}
[c]{c}%
A_{12}e_{n+1}\\
A_{22}e_{n+1}\\
A_{32}e_{n+1}\\
A_{42}e_{n+1}
\end{array}
\right)  \oplus \left(
\begin{array}
[c]{c}%
w_2(n)A_{12}e_{n}\\
w_2(n)A_{22}e_{n}\\
w_2(n)A_{32}e_{n}\\
w_2(n)A_{42}e_{n}
\end{array}
\right)  
\in E_3
$$
Therefore the rest of the proof is as same as 
\ref{lemma;triangular}. 
\end{proof}

\begin{lemma} Suppose that $A_{11} = O$. 
Then $A_{32} = O$ and $A_{22} = \lambda I + N$ for some 
$\lambda \in {\mathbb C}$ and  an upper 
triangular matrix $N \in B(\ell^{2}({\mathbb N}))$ 
with zero diagonal. 
\end{lemma} 
\begin{proof}Since $w_2$ appears in $Vu$ instead of $w_1$, 
a diagonal block $A_{22}$ plays a similar role of  a diagonal 
block $A_{11}$ in the argument of 
the proof in Lemma  \ref{lemma;diagonal}. The rest is 
similarly proved as the first column  of the operator 
matrix $A =(A_{ij})_{ij}$ is zero. 
\end{proof}

\begin{lemma}Suppose that $A_{11} = O$. Then 
$A_{22} = O$. 
\label{lemma;A12}
\end{lemma}
\begin{proof}
Since $A$ is an idempotent and $A_{32} = A_{42} = O$, 
$A_{22}$ is also an idempotent. 
Thus $A_{22} = O$ or $A_{22} = I$ by \ref{lemma;known}. 
It is enough to show  that $A_{22} \not= I$. 
On the contrary suppose that $A_{22} = I$. 
Then 
$$
V((0,e_1,0,0)\oplus (e_1,0,0,0)) 
= (A_{12}e_1,e_1,0,0)\oplus (0,0,0,0) \in E_3.
$$   
This implies that $A_{12}(21)=-\frac{1}{w_{1}(1)}.$
Since
$V((0,e_{n+1},0,0)\oplus(e_{n+1},w_2(n)e_n,0,0)) \in E_3$, 
we have 
\begin{align*}
& \left(
\begin{array}
[c]{c}%
A_{12}e_{n+1}\\
e_{n+1}\\
0\\
0
\end{array}
\right)  \oplus \left(
\begin{array}
[c]{c}%
w_2(n)A_{12}e_{n}\\
w_2(n)e_{n}\\
0\\
0
\end{array}
\right)  \\
& = 
\left(
\begin{array}
[c]{c}%
(x_{1}(k))_k\\
(x_{2}(k))_k\\
(x_{3}(k))_k\\
(x_{4}(k))_k
\end{array}
\right)  \oplus \left(
\begin{array}
[c]{c}%
(x_{1}(k+1)w_{1}(k)+x_{2}(k))_k\\
(x_{2}(k+1)w_{2}(k)+x_{3}(k))_k\\
(x_{3}(k+1)w_{3}(k)+x_{4}(k))_k\\
(y, (x_4(k))_k)
\end{array}
\right)  
\in E_3.
\end{align*}
for some 
$x_1,x_2,x_3,x_4 \in \ell^{2}({\mathbb N})$ and 
$y \in {\mathbb C}$. 

Then $A_{12}(n+2,n+1) = x_1(n+2)$ and $x_2(n+1) = 1$.
We also have 
$$
w_2(n)A_{12}(n+1,n) = x_1(n+2)w_1(n+1) + x_2(n+1).
$$
Therefore 
$$
A_{12}(n+2,n+1)
=\frac{w_{2}(n)}{w_{1}(n+1)} A_{12}(n+1,n)-\frac{1}{w_{1}(n+1)}.
$$
Hence we have 
$$
A_{12}(21)=-\frac{1}{w_{1}(1)}, \ \ \ 
A_{12}(32)=-\frac{w_{2}(1)}{w_{1}(1)w_{1}(2)}-\frac{1}{w_{1}(2)}, 
$$
$$
A_{12}(43)=-\frac{w_{2}(2)w_{2}(1)}{w_{1}(3)w_{1}(2)w_{1}(1)}
-\frac{w_{2}(2)}{w_{1}(3)w_{1}(2)} -\frac{1}{w_{1}(3)}, \ \dots 
$$
As $w_1(n) > 0$ and $w_2(n) > 0$, 
$$
|A_{12}(n+2,n+1)| \geq \frac{\prod_{k=1}^{n}w_2(k)}
{w_1(n+1)\prod_{k=1}^{n}w_1(k)}.
$$
Since $ 1 < w_1(n) \leq 4$ and 
 $\limsup_{n \rightarrow \infty} \prod_{k=1}^{n} 
\frac{w_2(k)}{w_1(k)} = \infty$ by Lemma \ref{lemma;sequences}, 
we have $\limsup_{n \rightarrow \infty}|A_{12}(n+2,n+1)|
= \infty $. This contradicts the fact that 
$||A_{12}|| < \infty. $
Therefore $A_{22} \not= I$. Hence  $A_{22} = O$.  
\end{proof}

\begin{lemma}Suppose that $A_{11} = O$. Then 
$A_{12} = O$, $A_{43} = O$, $A_{33} = O$, 
$A_{23} = O$ and $A_{13} = O.$
\end{lemma}
\begin{proof}Similar arguments before  show that 
$A_{12} = O$, $A_{43} = O$.  
$A_{33}$ is an idempotent and 
$A_{33} = \lambda I + N$ for some 
$\lambda \in {\mathbb C}$ and  an upper 
triangular matrix $N \in B(\ell^{2}({\mathbb N}))$ 
with zero diagonal. 
Thus $A_{33} = O$ or $A_{33} = I$ by Lemma \ref{lemma;known}. 
It is enough to show  that $A_{33} \not= I$. 
On the contrary suppose that $A_{33} = I$. 
Then 
$$
A_{23}(21)=-\frac{1}{w_{2}(1)}, \ \ \ 
A_{23}(n+2,n+1)
=\frac{w_{3}(n)}{w_{2}(n+1)} A_{23}(n+1,n)-\frac{1}{w_{2}(n+1)}.
$$
As in the proof of Lemma \ref{lemma;A12}, we have  
$\limsup_{n \rightarrow \infty}|A_{23}(n+2,n+1)|
= \infty $. This contradicts the fact that 
$||A_{23}|| < \infty. $
Therefore $A_{33} \not= I$. Hence  $A_{33} = O$. 
The rest is similarly proved. 
\end{proof}

\begin{lemma}Suppose that $A_{11} = O$. Then 
$A_{44} = O$, $A_{34} = O$, $A_{24} = O$ 
and $A_{14} = O$.
\end{lemma}
\begin{proof}Since the fourth column of operator matrix 
$T = (T_{ij})_{ij}$ has a different form than the the other columns, 
we need to be careful to investigate.  

Since $u = (0,0,0,0)\oplus(0,0,0,e_1) \in E_3$, we have 
\begin{align*}
Vu & =  
\left(
\begin{array}
[c]{c}%
0\\
0\\
0\\
0
\end{array}
\right) \oplus \left(
\begin{array}
[c]{c}%
A_{14}e_1\\
A_{24}e_1\\
A_{34}e_1\\
A_{44}e_1
\end{array}
\right)  \\
& = 
\left(
\begin{array}
[c]{c}%
(x_{1}(k))_k\\
(x_{2}(k))_k\\
(x_{3}(k))_k\\
(x_{4}(k))_k
\end{array}
\right)  \oplus \left(
\begin{array}
[c]{c}%
(x_{1}(k+1)w_{1}(k)+x_{2}(k))_k\\
(x_{2}(k+1)w_{2}(k)+x_{3}(k))_k\\
(x_{3}(k+1)w_{3}(k)+x_{4}(k))_k\\
(y, (x_4(k))_k)
\end{array}
\right)  
\in  E_3.
\end{align*}
for some 
$x_1,x_2,x_3,x_4 \in \ell^{2}({\mathbb N})$ and 
$y \in {\mathbb C}$. Then $x_1 = x_2 = x_3 = x_4 = 0$.  
Therefore  $A_{14}(k,1) = A_{24}(k,1) = A_{34}(k,1) = 0$ 
for any $k \geq 1$. 
We also have $A_{44}(k,1) = 0$ for any $k \geq 2$. 

Since $u = (0,0,0,e_n)\oplus(0,0,e_n,e_{n+1}) \in E_3$, 
we have 
\begin{align*}
Vu & =  
\left(
\begin{array}
[c]{c}%
A_{14}e_n\\
A_{24}e_n\\
A_{34}e_n\\
A_{44}e_n
\end{array}
\right)  \oplus \left(
\begin{array}
[c]{c}%
A_{14}e_{n+1}\\
A_{24}e_{n+1}\\
A_{34}e_{n+1}\\
A_{44}e_{n+1}
\end{array}
\right)  \\
& = 
\left(
\begin{array}
[c]{c}%
(x_{1}(k))_k\\
(x_{2}(k))_k\\
(x_{3}(k))_k\\
(x_{4}(k))_k
\end{array}
\right)  \oplus \left(
\begin{array}
[c]{c}%
(x_{1}(k+1)w_{1}(k)+x_{2}(k))_k\\
(x_{2}(k+1)w_{2}(k)+x_{3}(k)_k\\
(x_{3}(k+1)w_{3}(k)+x_{4}(k))_k\\
(y, (x_4(k))_k))
\end{array}
\right)  
\in E_3.
\end{align*}
for some 
$x_1,x_2,x_3,x_4 \in \ell^{2}({\mathbb N})$ and 
$y \in {\mathbb C}$. 

Then 
$$
A_{44}(k+1,n+1) = x_4(k) = A_{44}(k,n) \text{ for any } 
k \geq 1, \ n \geq 1. 
$$
Since $A_{44}(1,1) = y$ and $A_{44}(k,1) = 0$ for $k \geq 2$, 
$A_{44} = yI + N$ for some $y \in {\mathbb C}$ and 
an upper triangular matrix $N$ 
with zero diagonal.  Since $A_{44}$ is an idempotent, 
$A_{44} = O$ or $A_{44} = I$. We shall show that 
$A_{44} \not= I$. On the contrary assume that 
$A_{44} = I$. Then 
$x_4 = A_{44}e_n = e_n$.  Moreover 
$$
A_{34}(k,n+1) = x_3(k+1)w_3(k) + e_n(k) = A_{34}(k+1,n)w_3(k) + e_n(k).
$$
This implies that 
$$
A_{34}(e_{n+1}) = B_{w_3}A_{34}(e_n) + e_n .
$$
Since $A_{34}(e_1) = 0$, we have 
$A_{34}(e_2) = e_1$, $A_{34}(e_3) = e_2$, \\
$A_{34}(e_4) = w_3(1)e_1 + e_3$, 
$A_{34}(e_5) = w_3(2)e_2 + e_4$, \\
$A_{34}(e_6) = w_3(1)w_3(2)e_1 + w_3(3)e_3 + e_5$, \dots .
Therefore 
$$
A_{34} =\left(
\begin{array}
[c]{cccccccccc}%
0 & 1 & 0 & w_3(1) & 0      & w_3(1)w_3(2) & 0  & w_3(1)w_3(2)w_3(3) & \dots \\
0 & 0 & 1 & 0      & w_3(2) & 0            & w_3(2)w_3(3) &0 & \dots\ \\
0 & 0 & 0 & 1      & 0      & w_3(3)       & 0 & w_3(3)w_3(4) & \dots \\
0 & 0 & 0 & 0      & 1      & 0            & w_3(4)       & 0 & \dots \\
0 & 0 & 0 & 0      & 0  &  1      & 0      & w_3(5)       & \dots \\
\vdots & \vdots & \vdots & \vdots & \vdots & \vdots &\vdots 
&\vdots &\vdots 
\end{array}
\right)
$$
In particular, we have 
$$
A_{34}(1,2n) = \prod_{k=1}^{n-1}w_3(k).
$$
Since 
$\lim_{n \rightarrow \infty} \prod_{k=1}^{n} w_3(k) = \infty$ ,
$\lim_{n \rightarrow \infty} A_{34}(1,2n) = \infty$.  
This contradicts that $A_{34}$ is bounded. 
Therefore $A_{44} = O$. 
Moreover 
\begin{align*}
A_{34}(k,n+1) & = x_3(k+1)w_3(k) + x_4(k) \\
  & = x_3(k+1)w_3(k) + A_{44}e_n = A_{34}(k+1,n)w_3(k).
\end{align*}
Since $A_{34}(k,1) = 0$, we have $A_{34} = O$. 
Similarly we have $A_{24} = O$ and $A_{14} = O. $
\end{proof}

The proof of Theorem \ref{thm:exotic examples} will be 
completed by the following Lemma: 

\begin{lemma}The system $\mathcal{S}_{w,N}$ of four subspaces 
is indecomposable. 
\end{lemma}
\begin{proof}Let $V \in End (\mathcal{S}_{w,N})$ satisfy $V^2 = V$ as 
in the beginning of the proof of Theorem \ref{thm:exotic examples}. 
Then $V = A \oplus A$ and $A_{11} = O$ or $A_{11} = I$ 
by  Lemma \ref{lemma;diagonal}. If $A_{11}= O$, then $A = O$ 
by the preceding 
lemmas so that $V = O$.  If $A_{11}$ = I,  then 
$(I-A)_{11} = 0$. By replacing 
$V$ by an idempotent $I-V$, the same argument implies 
$I-A = O$, so that $V = I$. 
This establishes that  $\mathcal{S}_{w,N}$ is indecomposable. 
\end{proof}

\section{being exotic }

In this section we shall show that the indecomposable systems 
$\mathcal{S}_{w,N}$ constructed in the preceding section are exotic 
in the sense that $\mathcal{S}_{w,N}$ are not 
isomorphic to any closed operator system ${\mathcal S}_A$ under any 
permutation of subspaces. 
We recall a necessary criterion in \cite{EW}. 

\bigskip
\noindent
{\bf Definition}(intersection diagram). Let 
$\mathcal{S} =(H;E_1,E_2,E_3,E_4)$ be a system of four
subspaces. The {\it intersection diagram} for a system $\mathcal{S}$
is an undirected graph 
$\Gamma_{\mathcal{S}} 
= (\Gamma_{\mathcal{S}}^0,\Gamma_{\mathcal{S}}^1)$ 
with the set of vertices $\Gamma_{\mathcal{S}}^0$ and the 
set of edges $\Gamma_{\mathcal{S}}^1$ defined by 
$\Gamma_{\mathcal{S}}^0 = \{1,2,3,4\}$ 
and for $i \not= j \in \{1,2,3,4\}$ 
\[
\circ _i  \ ^{\line(1,0){20}}  \circ _j \ \text{ if and only if } 
E_i \cap E_j = 0 .
\]

\begin{lemma}(\cite[Lemma 10.4]{EW})
Let 
$\mathcal{S} = \mathcal{S}_{T,S} =(H;E_1,E_2,E_3,E_4)$ be a closed 
operator system. Then the intersection diagram $\Gamma_{\mathcal{S}}$ 
for the system $\mathcal{S}$ contains 
\[
\circ _4 \ ^{\line(1,0){20}} \circ _1 \ ^{\line(1,0){20}} 
\circ _2 \ ^{\line(1,0){20}} \circ _3 \ , 
\]
that is, $E_4 \cap E_1 = 0$, $E_1 \cap E_2 = 0$ and 
$E_2 \cap E_3 = 0$.  In particular, then the intersection diagram 
$\Gamma_{\mathcal{S}}$ is a connected graph.  
\end{lemma}

\begin{prop}
The indecomposable systems 
$\mathcal{S}_{w,N}$ constructed in the preceding section are 
not isomorphic to any closed operator systems under any 
permutation of subspaces. 
\end{prop}
\begin{proof}
It is clear that $E_4 \cap E_1 = 0$, $E_1 \cap E_2 = 0$ and 
$E_2 \cap E_4 = 0$.  
Since $(e_1,0,...,0)\oplus(0,0,...,0) \in E_1 \cap E_3$, 
we have $E_1 \cap E_3 \not= 0$.  Because 
$(0,0,...,0)\oplus(0,0,...,e_1)  \in E_2 \cap E_3$, 
we have $E_2 \cap E_3 \not= 0$. By Lemma \ref{lemma;sequences}, 
there exists a non-zero vector $x_1 \in K = \ell ^2(\mathbb N)$ 
with $B_{w_1}x_1 = x_1$. 
Then $(x_1,0,...,0,0)\oplus(x_1,0,...,0,0) \in E_3 \cap E_4$, 
so that $E_3 \cap E_4 \not= 0$.  Therefore the vertex $3$ is 
not connected to any other vertices $1, 2, 4$.  Thus 
the intersection diagram 
$\Gamma_{\mathcal{S}_{w,N}}$ is not a connected graph. This implies 
that $\mathcal{S}_{w,N}$ 
is not isomorphic to any closed operator system under any 
permutation of subspaces. 
\end{proof}

\section{defect computation}

We shall compute the defect of the indecomposable systems 
$\mathcal{S}_{w,N}$ constructed in section 3. 

\begin{lemma}
For fixed $j \in {\mathbb N}$ and $b \in \ell^2({\mathbb N})$, consider 
an equation $B_{w_j} u  + b = u$ for
unknown sequence $u \in \ell^2({\mathbb N})$.
Suppose that there exists a polynomial $p(t)$ of degree $r$ 
with positive coefficients such that 
$|b(n+1)| \leq p(n)(\frac{3}{4})^n$ for  $n \in {\mathbb N}$. 
For any $c \in {\mathbb C}$, put  $u(1) = c$ and let 
$$
u(n+1) = \frac{c}{\prod_{k=1}^{n}w_j(k)}
-\sum_{m=1}^{n}\frac{b(m)}{\prod_{k=m}^{n}w_j(k)} 
\text{ for } n \in {\mathbb N}.
$$
Then there exists a polynomial $q(t)$ of degree $r+1$ such that 
$u := (u(n))_n$ satisfies 
$|u(n+1)| \leq q(n)(\frac{3}{4})^n$ for any $n$. 
Moreover  $u$ is in $\ell^2({\mathbb N})$ and 
a solution of the equation $B_{w_j} u  + b = u$. Conversely
any solution $u$ has this form. 
\end{lemma}
\begin{proof}
The  equation $B_{w_j} u  + b = u$ implies that 
$$
w_j(n)u(n+1) + b(n) = u(n)  \text{ for } n \in {\mathbb N}.
$$
Hence $u(n+1) = \frac{u(n)}{w_j(n)}-\frac{b(n)}{w_j(n)}$. 
Therefore any solution $u$ has the desired form:
$$
u(n+1) = \frac{c}{\prod_{k=1}^{n}w_j(k)}
-\sum_{m=1}^{n}\frac{b(m)}{\prod_{k=m}^{n}w_j(k)} 
\text{ for } n \in {\mathbb N} .
$$

Since  $\frac{4}{3} \leq w_j(n) \leq 4$  
$$
|u(n+1)| 
\leq |c|(\frac{3}{4})^n + \sum_{m=1}^{n}p(m-1)(\frac{3}{4})^{m-1}
(\frac{3}{4})^{n-m+1} \leq q(n)(\frac{3}{4})^n
$$
for some  polynomial $q(t)$ of degree $r+1$. 
It is easy to see that $u$ satisfies the equation and is in 
$\ell^2({\mathbb N}).$
\label{lemma;equation}
\end{proof}

\begin{prop} For any natural number $N$  
the indecomposable systems $\mathcal{S}_{w,N}$ have the defect 
$\rho(\mathcal{S}_{w,N}) = \frac{2N+1}{3}.$
\end{prop}
\begin{proof}
We need to compute  $\dim (E_{i}\cap E_{j})$ and 
$\dim ((E_i + E_j)^{\perp}).$ 
It is obvious that $\dim (E_i\cap E_j) = 0$ and 
$\dim ((E_i + E_j)^{\perp}) = 0 $ for any $i, j = 1,2,4$ with 
$i \not=j$.  

We consider $E_2 + E_3$. Since $E_2 = 0 \oplus K$ and 
$E_3 \supset \{(x\oplus Tx) ; x \in K\}$, 
$E_2 + E_3 \supset H$. Thus $\dim ((E_2 + E_3)^{\perp}) = 0$. 
Next we investigate $E_2 \cap E_3$. We see that 
$\{(0,0,...,0,0)\oplus(0,0,...,0,\alpha e_{1}) 
; \alpha \in {\mathbb C} \} \subset E_{2}\cap E_{3}$. 
Conversely take any 
$$
\left(
\begin{array}
[c]{c}%
0\\
\vdots\\
0\\
0
\end{array}
\right)  \oplus\left(
\begin{array}
[c]{c}%
z_{1}\\
\vdots\\
z_{N}\\
z_{N+1}%
\end{array}
\right)  
= 
\left(
\begin{array}
[c]{c}%
x_{1}\\
\vdots\\
x_{N}\\
y
\end{array}
\right)  \oplus\left(
\begin{array}
[c]{c}%
B_{w_{1}}x_{1}+x_{2}\\
\vdots\\
B_{w_{N}}x_{N}+y\\
Sy+\alpha e_{1}%
\end{array}
\right) 
\in E_{2}\cap E_{3}.
$$
Since $y=0,x_{1}=x_{2}=\cdots=x_{N}=0,$  $z_1= \dots = 0$ and 
$$
E_{2}\cap E_{3}=\{(0,0,...,0,0)\oplus(0,0,...,0,\alpha e_{1}) 
; \alpha \in {\mathbb C} \}.
$$
Therefore  $\dim (E_{2}\cap E_{3}) =1.$ 
Next we shall show that $E_{1}+E_{3}=H.$ Since 
$E_3 \supset \graph T$, $E_1 + E_3 \supset 0 \oplus \Im T$. 
And $\Im T = (L,L,\dots, L,\Im S)$, because $B_{w_k}$ is onto. 
Considering one dimensional perturbation by $\{(0,0,...,0,0)\oplus(0,0,...,0,\alpha e_{1}) 
; \alpha \in {\mathbb C} \}$, $E_{1}+E_{3}= H$, so that 
$\dim ((E_{1}+E_{3})^{\perp}) = 0$. 
Consider $E_1 \cap E_3$.  
Take any 
$$
\left(
\begin{array}
[c]{c}%
x_{1}\\
\vdots\\
x_{N}\\
y
\end{array}
\right)  \oplus\left(
\begin{array}
[c]{c}%
0\\
\vdots\\
0\\
0
\end{array}
\right)  =\left(
\begin{array}
[c]{c}%
x_{1}\\
\vdots\\
x_{N}\\
y
\end{array}
\right)  \oplus\left(
\begin{array}
[c]{c}%
B_{w_{1}}x_{1}+x_{2}\\
\vdots\\
B_{w_{N}}x_{N}+y\\
Sy+\alpha e_{1}
\end{array}
\right)  \in E_{1}\cap E_{3}.
$$

Then $y=0,\alpha=0.$ Since $B_{w_{N}}x_N = 0$, 
$x_{N}=(x_{N}(1),0,0,0,...).$ 
From  $B_{w_{N-1}}x_{N-1} + x_N = 0$, we have 
$x_{N-1}=(x_{N-1}(1),-\frac{x_{N}(1)}{w_{N-1}(1)},0,0,0,...).$ \\
We continue in this way to obtain 
$x_{N-2}=(x_{N-2}(1),-\frac{x_{N-1}(1)}{w_{N-2}(1)},\frac{x_{N}(1)}
{w_{N-2}(2)w_{N-1}(1)},0,0,...),$

$x_{N-3}=(x_{N-3}(1),-\frac{x_{N-2}(1)}{w_{N-3}(1)},\frac{x_{N-1}(1)}
{w_{N-3}(2)w_{N-2}(1)},-\frac{x_{N}(1)}{w_{N-3}(3)w_{N-2}(2)w_{N-1}
(1)},0,0,...),$

......, and 

$x_{1}=(x_{1}(1),(-1)^{1}\frac{x_{2}(1)}{w_{1}(1)},(-1)^{2}\frac{x_{3}%
(1)}{w_{1}(2)w_{2}(1)},(-1)^{3}\frac{x_{4}(1)}{w_{1}(3)w_{2}(2)w_{3}(1)},...,$

$(-1)^{N-1}\frac{x_{N}(1)}{w_{1}(N-1)w_{2}(N-2)\cdots w_{N-1}(1)},0,0,...)$

Conversely for any parameters 
$x_{1}(1),x_{2}(1),...,x_{N}(1) \in {\mathbb C}$, vectors 
$x_{1},x_{2},...,x_{N}$  with the above forms and 
$y = 0, \alpha = 0$ give elements in $E_1 \cap E_3$. 
Therefore  $\dim (E_{1}\cap E_{3})=N.$

Next we investigate $(E_3 + E_4)^{\perp}$.
Since $E_4^{\perp} = \{(-y,y) \in H ; y \in K\}$ and 
$(\graph T)^{\perp} = \{(-T^*z,z) \in H ; z \in K\}$, 
we have 
$$
E_3^{\perp} \cap E_4^{\perp}
= \{(-T^*z,z) \in H ; z \in K, \ T^*z = z, (z|(0,...,0,e_1)) = 0 \}, 
$$
Let $z = (z_1,\dots, z_N,w) \in K$. 
Since $B_{w_{k}}^* = S_{\overline{w_k}} = S_{{w_k}}$ is a 
weighted shift,  $T^*z = z$ implies that 
$$
(S_{w_1}z_1,z_1 + S_{w_2}z_2, \dots ,z_{N-1} + S_{w_N}z_N, 
z_N + S^*w) = (z_1,\dots, z_N,w).
$$
From $S_{w_1}z_1 = z_1$, we have $z_1 = 0$. 
Then $S_{w_2}z_2 = z_2$. Hence $z_2 = 0$. We continue in this way 
to obtain $z_3 = \dots = z_N = 0$. Therefore 
$0 = (z|(0,...,0,e_1)) = (w|e_1)$. Furthermore $S^* w = w$. 
Hence $w = 0$. Thus $z = 0$. Hence 
$E_3^{\perp} \cap E_4^{\perp} = 0$. 

Finally we investigate $E_{3}\cap E_{4}$. 
Take any 
$$
\left(
\begin{array}
[c]{c}%
x_{1}\\
\vdots\\
x_{N}\\
y
\end{array}
\right)  \oplus\left(
\begin{array}
[c]{c}
B_{w_{1}}x_{1}+x_{2}\\
\vdots\\
B_{w_{N}}x_{N}+y\\
Sy+\alpha e_{1}%
\end{array}
\right)
=  
\left(
\begin{array}
[c]{c}%
x_{1}\\
\vdots\\
x_{N}\\
y
\end{array}
\right)  \oplus\left(
\begin{array}
[c]{c}%
x_{1}\\
\vdots\\
x_{N}\\
y
\end{array}
\right)
\in E_{3} \cap E_{4}.
$$
Since $Sy + \alpha e_1 = y$, $y = (\alpha,\alpha,\alpha,\dots ),$  
As $y \in \ell^2({\mathbb N})$, $\alpha = 0$ and $y = 0$. 
Then $B_{w_N}x_N = x_N$. Hence $w_N(n)x_N(n+1) = x_N(n)$ for 
$n \in {\mathbb N}$. Therefore there exists a constant $c_N$ 
such that  
$$
x_{N}=c_{N}(1,\frac{1}{w_{N}(1)},\frac{1}{w_{N}(2)w_{N}(1)},
\frac{1}{w_{N}(3)w_{N}(2)w_{N}(1)},\cdots ).
$$  
Then $|x_N(n+1)| = |c_N| \frac{1}{\prod_{k=1}^{n}w_N(k)}
\leq |c_N|(\frac{3}{4})^n$. Thus $x_N \in \ell^2({\mathbb N})$. 
Apply Lemma \ref{lemma;equation} for 
the equations $B_{w_{j}}x_{j}+x_{j+1} = x_j$ for 
$j = N-1,N-2, \dots, 1$ step by step. 
There exist parameters $c_{N-1}, c_{N-2}, \dots, c_1$ such 
that $x_j(1) = c_j$  for 
$j = N,N-1, \dots, 1$ and the other components $x_j(n)$ for 
$n \geq 2$ are uniquely determined by these parameters. In fact,
$$
x_{j}(n+1)  = \frac{c_j}{\prod_{k=1}^{n}w_{j}(k)}
-\sum_{m=1}^{n}\frac{x_{j+1}(m)}{\prod_{k=m}^{n}w_{j}(k)} 
\text{ for } n \in {\mathbb N} .
$$
Conversely any $x_1, \dots, x_N$ with this form gives an element of 
$E_{3}\cap E_{4}$. 
Hence $\dim (E_{3}\cap E_{4}) = N$.  
Therefore $\rho(\mathcal{S}_{w,N}) = \frac{2N+1}{3}.$
\end{proof}

Finally we get the following theorem.

\begin{thm}
There exist  exotic indecomposable systems of 
four subspaces with the defect $\frac{2n+1}{3}(n\in {\mathbb N}).$
\end{thm}

\end{document}